\documentclass[11pt]{article}
\usepackage[cp1251]{inputenc}
\usepackage[russian, english]{babel}
\usepackage{amsfonts,amsmath,enumerate}
\usepackage{graphicx,epsfig}

\newtheorem{theorem}{Theorem}

\def\beq{\begin{equation}}
\def\eeq{\end{equation}}
\def\R{{\mathbb R}}

\def\dn{\mathrm{dn}\,}
\def\sn{\mathrm{sn}\,}
\def\cn{\mathrm{cn}\,}
\def\am{\mathrm{am}\,}
\def\E{\mathrm{E}}
\def\solv{{\mathrm{SOLV^-}}}
\def\solve{{\mathrm{SOLV^+}}}
\def\sol{{\mathrm{Sol}}}
\def\se{{\mathrm{SE(2)}}}
\def\sh{{\mathrm{SH(2)}}}

\begin{document}
\title{Sub-Riemannian geodesics on the three-dimensional solvable non-nilpotent Lie group $\solv$}
\author{Mazhitova A.D.
\thanks{
Mechanical-Mathematical Department, Al-Farabi Kazakh National
University, Al-Farabi ave. 71, Almaty 050038, Kazakhstan; e-mail:
Akmaral.Mazhitova@kaznu.kz} }
\date{}
\maketitle

\section{Introduction}

In this paper we describe the geodesics of a left-invariant
sub-Riemannian metric on the three-dimensional solvable Lie group
$\mathrm{Sol}$.

This group is widely known in geometry, because it allows compact
quotient-spaces and it gives one of the Thurston three-dimensional
geometries \cite{Thurston}. By the classification theorem of
Agrachev--Barilari \cite{AB} there are invariant sub-Riemannian
geometries realized on four solvable non-nilpotent Lie groups:
$\se$, $\sh$, $\solv$, and $\solve$.

In this classification, our geometry corresponds to the case
$\solv$:
$$
\sol = \solv.
$$
The case of $\mathrm{SOLV}^+$ we shall consider separately.

Various aspects of the integration of geodesic flows on
sub-Riemannian manifolds has been widely studied (see, for example,
\cite{MS,S,AS,Sachkov,T}). Note, that the geodesics of other
three-dimensional nonsolvable or nilpotent sub-Riemannian geometries
have been described recently in terms of elementary functions
\cite{BR,CM}. In our situation it is necessary to use elliptic
functions.

We thank I.A. Taimanov for posing the problem and Ya.V. Bazaikin for
helpful discussions.

\section{Basic definitions}
\subsection{Geodesics of sub-Riemannian manifolds}

Let $M^{n}$ be a smooth $n$-dimensional manifold. A smooth family of
$k$-dimensional subspaces in the tangent spaces to points of $M^{n}$
$$
\Delta=\left\{\Delta(q):\ \ \Delta(q) \in T_{q}M^{n} \ \ \ \
\forall q \in M^{n},\ \ \dim \Delta(q)=k\right\}
$$
is called completely nonintegrable, if the vector fields tangent to
$\Delta$, and all their iterated commutators generate the tangent
bundle $TM^n$:
$$
\mathrm{span} \left\{[f_{1},[\dots[f_{m-1},f_{m}]\dots]](q): f_{i}(p) \in
\Delta(p)\ \ \forall p \in M^{n}, m=1,\dots \right\}=T_{q}M^{n}.
$$
Sometimes this distribution is called completely nonholonomic.

A two-dimensional distributions in a three-dimensional manifold is
completely nonholonomic if and only if
$$
\mathrm{span}\{f_{1}(q), f_{2}(q), [f_{1}(q), f_{2}(q)]\}=T_{q}M^3,
$$
where at every point $q$ the vectors $f_1(q)$ and $f_2(q)$ form a
basis in $\Delta(q)$.

Let $g_{ij}$ be a complete Riemannian metric on $M^n$. A triple
$(M^{n},\ \ \Delta, \ \ g_{ij})$ is called a sub-Riemannian
manifold. A Lipschitz continuous curve $\gamma:[0,T]\rightarrow
M^{n}$ is called admissible if $\dot{\gamma}(t)\in
\Delta(\gamma(t))$  for almost all $t\in [0,T]$. The length of this
curve is equal to
$$
l(\gamma)=\int_{0}^{T}\sqrt{g_{\gamma(t)} \left( \dot{\gamma}(t),
\dot{\gamma}(t) \right) }dt.
$$
The distance between two points on the manifold is defined by
formula
$$
d(q_{0},q_{1})=\inf_{\gamma\in\Omega_{q_{0},q_{1}}}l(\gamma),
$$
where $\Omega_{q_{0},q_{1}}$ is the set of all admissible curves
connecting points $q_{0}$ and $q_{1}$. This function
$d(\cdot,\cdot)$ is called the sub-Riemannian metric on $M^{n}$. A
geodesic of this metric is an admissible curve
$\gamma:[0,T]\rightarrow M^{n}$, which locally minimizes the length
functional $l(\gamma)$.

Geodesics of sub-Riemannian metrics satisfy the Pontryagin maximum
principle (see, for instance, \cite{AS}),
 which we formulate below.
Let $f_1,\dots,f_k$ be vector fields which are tangent to $\Delta$
and span $\Delta$ at every point of $M^n$ (or of a domain of $M^n$).

{\sc The Pontryagin maximum principle} is as follows:

\begin{itemize}
\item
{ \sl Let $M^{n}$ be a smooth $n$-dimensional manifold. Let us
consider for Lipschitz continuous curves the following minimum
problem
$$
\dot{q}=\sum_{i=1}^{k} u_{i}f_{i}(q),\ \ u_{i} \in \R,\ \
\int_{0}^{T}\sum_{i=1}^{k} u_{i}^{2}(t)dt\longrightarrow \min, \ \
q(0)=q_{0},\ \ q(T)=q_{1}
$$
with a fixed $T$. Let us consider the mapping $\mathcal{H}:T^{*}M^{n}\times \R\times
\R^{k} \rightarrow \R$, given by the function
$$\mathcal{H}(q, \lambda, p_{0}, u):= \langle \lambda,
\sum_{i=1}^{k}u_{i}f_{i}(q) \rangle + p_{0}\sum_{i=1}^{k}
u_{i}^{2}.
$$
If a curve $q(\cdot):[0,T]\rightarrow M^n$ with a control
$u(\cdot):[0,T]\rightarrow \R^{k}$ is optimal, then there exists Lipshitzian covector function
$\lambda(\cdot):t\in[0,T]\mapsto \lambda (t)\in T_{q(t)}^{*}M^{n}$, $(\lambda (t), p_0)\neq 0$
and a constant  $p_{0}\leq 0$  such that

i) $\dot{q}(t)=\dfrac{\partial \mathcal{H}}{\partial \lambda}(q(t),
\lambda(t), p_{0}, u(t)),$

ii) $\dot{\lambda}(t)=-\dfrac{\partial \mathcal{H}}{\partial q}(q(t),
\lambda(t), p_{0}, u(t)),$

iii) $\dfrac{\partial \mathcal{H}}{\partial u}(q(t),
\lambda(t), p_{0}, u(t))=0.$
}
\end{itemize}

A curve $q(\cdot):[0,T]\rightarrow M^{n}$, satisfying  the
Pontryagin maximum principle is called an extremal (curve). To such
a curve there corresponds a set of pairs $(\lambda(\cdot),p_{0})$.
The type (normal or abnormal) of an extremal depends on the value of
$p_0$:
\begin{itemize}
\item
if $p_{0}\neq0$, then the extremal is called \emph{normal};

\item
if $p_{0}=0$, then the extremal is called \emph{abnormal};
\item
extremal is called \emph{strictly abnormal} if it is not projected
(on $M^n$) onto a normal extremal.
\end{itemize}

For a normal extremal we may put $p_{0}=-\frac{1}{2}$.

\emph{Normal extremals} are geodesics \cite{AS}. In the contact
case, when at every point the distribution $\Delta$ coincides with
the annihilator of the contact form on $M^n$, there are no nontrivial
 abnormal extremals (this fact is indicated in \cite{BR}). In the case, when
the space of vector fields on a manifold is generated by vector
fields tangent to the nonholonomic distribution and their
commutators, there are no strictly abnormal extremals \cite{AS}. We
recall, that an extremal is called strictly abnormal, if its
projection on $M^n$ does not coincide with the projection of any
normal extremal. Both of the above statements apply to
three-dimensional sub-Riemannian manifold $M^3$.

By iii), $u_{i}= \langle \lambda(t),f_{i}(t) \rangle$ and a curve
$q(\cdot):[0,T]\rightarrow M^{n}$ is geodesic if and only if it is
the projection onto $M^{n}$ of a solution $(\lambda(t),q(t))$ of a
Hamiltonian system on $T^{*}M^{n}$ with the following Hamiltonian
function:
$$
H(\lambda,q)=\frac{1}{2}\left(\sum_{i=1}^{k} \langle \lambda, F_{i}
\rangle^{2}\right),\ \ q\in M^{n},\ \ \lambda \in T_{q}^{*}M^{n}.
$$
The Hamiltonian $H$ is constant along any
solution of the Hamiltonian system. Moreover, $H=\frac{1}{2}$
if and only if the geodesic is length parameterized.

\subsection{Elliptic functions. Jacobi functions}

Let us recall some necessary facts of Jacobi elliptic functions.
The integrals
$$
\int_{0}^{x} \frac{dx}{\sqrt{(1-x^{2})(1-k^{2}x^{2})}},
$$
and
$$
\int_{0}^{x} \frac{\sqrt{1-k^{2}x^{2}}}{\sqrt{(1-x^{2})}}dx,
$$
are called elliptic integrals of first and second kind,
respectively, in the normal Legendre form (see \cite{GR,Aks}), where
$k$ ($0<k<1$) is the modulus of these integrals, $k'=\sqrt{1-k^{2}}$
is the additional modulus. By the substitution $x=\sin \varphi$
these integrals reduce to the normal trigonometric form

\beq
\label{ell}
F(\varphi, k)=\int_{0}^{\varphi} \frac{d\alpha}{\sqrt{(1-k^{2}\sin^{2}\alpha)}}
=\int_{0}^{\sin \varphi} \frac{dx}{\sqrt{(1-x^{2})(1-k^{2}x^{2})}},
\eeq

\beq
\label{ell2}
\E(\varphi, k)=\int_{0}^{\varphi}\sqrt{1-k^{2}\sin^{2}\alpha} \ \ d \alpha
=\int_{0}^{\sin \varphi} \frac{\sqrt{1-k^2 x^2}}{\sqrt{1-x^2}} dx.
\eeq

Consider an integral of first kind in the normal trigonometric form

$$
v=\int_{0}^{\varphi}\frac{d\varphi}{\sqrt{1-k^{2}\sin^{2}\varphi}}.
$$
Now consider the upper limit as a function of $v$. This function is
denoted by
$$
\varphi=\am(v, k)=\am v
$$
and is called the amplitude, and this process is called an inversion
of the integral. Thus, next functions:
$$
\sin \varphi=\sin (\am v)=\sn v,
$$
$$
\cos \varphi=\cos (\am v)=\cn v,
$$
$$
\Delta  \am v=\sqrt{1- k^{2}\sin^{2} \varphi}=\sqrt{1-k^{2}\sn^{2} v}=\dn v
$$
are called the Jacobi functions and are related by
$$
\sn^{2}v+\cn^{2}v=1,\ \  \dn^{2}v+k^{2}\sn^{2}v=1.
$$
By derivation, we obtain
$$
\frac{d \sn v}{dv}=\cn v \dn v,
$$
$$
\frac{d \cn v}{dv}=-\sn v \dn v,
$$
$$
\frac{d \dn v}{dv}=-k^{2}\sn v \cn v
$$
and conclude, that
$$
\left( \frac{d \sn v}{dv}\right)^{2}=(1-\sn^{2}v)(1-k^{2}\sn^{2}v),
$$
\beq
\label{ell3}
\left(\frac{d \cn v}{dv}\right)^{2}=(1-\cn^{2}v)(k'^{2}+k^{2}\cn^{2}v),
\eeq
$$
\left(\frac{d \dn v}{dv}\right)^{2}=(1-\dn^{2}v)(\dn^{2}v-k'^{2}).
$$
The first equation of (\ref{ell3}) implies that $\sn v$ is
the inversion of the first kind elliptic integral in
the normal Legendre form
\beq
\label{ell4}
v=\int_{0}^{\sn v}\frac{dx}{\sqrt{(1-x^{2})(1-k^{2}x^{2})}}.
\eeq

From the second and third equations obtained that
$\cn v$ and $\dn v$ are the result
of conversion of next functions
\beq
\label{ell5}
v=\int_{1}^{\cn v}\frac{dx}{\sqrt{(1-x^{2})(k'^{2}+k^{2}x^{2})}},
\eeq
\beq
\label{ell6}
v=\int_{1}^{\dn v}\frac{dx}{\sqrt{(1-x^{2})(x^{2}-k'^{2})}}.
\eeq

All Jacobi functions are periodic.
Note that the function $\sn v$ is odd,
but $\cn v$ and $\dn v$ are even, therefore we assume, what
in the two last integral, when the functions $\cn v$ and $\dn v$
pass through the critical points,
respectively changes a sign of radical.

\section{Sub-Riemannian problem on the group $\solv$}

Let us consider the three-dimensional Lie group $\solv$ formed by
all matrices of the form
$$
\left(
\begin{array}{ccc}
 e^{-z}& 0 & x \\
0 & e^z & y \\
0 & 0 & 1
\end{array}
\right), x,y,z \in \R.
$$
Its Lie algebra is spanned by the vectors
$$
e_{1} = \left(
\begin{array}{ccc}
0 & 0 & 1 \\
0 & 0 & 0 \\
0 & 0 & 0
\end{array}
\right) ,\ \
e_{2} = \left(
\begin{array}{ccc}
0 & 0 & 0 \\
0 & 0 & 1 \\
0 & 0 & 0
\end{array}
\right),\ \
e_{3} = \left(
\begin{array}{ccc}
-1 & 0 & 0 \\
0 & 1 & 0 \\
0 & 0 & 0
\end{array}
\right),
$$
meeting the following commutation relations:
$$
[e_{1},e_{2}]=0; \ \ [e_{1},e_{3}]=e_{1}; \ \ [e_{2},e_{3}]=-e_{2}.
$$
We take a new basis
\beq
\label{e3}
a_{1}=e_{1}+e_{2}; \ \  a_{2}=e_{1}-e_{2}; \ \  a_{3}=e_{3},
\eeq
in which
the commutation relations take the form
$$
[a_{1},a_{2}]=0, \ \ [a_{1},a_{3}]=a_{2}, \ \  [a_{2},a_{3}]=a_{1}.
$$
Let us consider a left-invariant metric on $\solv$,
which is defined by its values in the unit of the group:
$$
(e_{i},e_{j})=\delta_{ij}.
$$

The Lie group $\solv$ is diffeomorphic to the space
$\R^3$. Indeed, $x, y, z$ are the global coordinates on $\solv$
and they also may be considered as global coordinates on
$\R^3$.
The tangent space at each point of $\solv$ is spanned by
matrices of the form
$$
\partial_{x} = \left(
\begin{array}{ccc}
0 & 0 &1 \\
0 & 0 & 0 \\
0 & 0 & 0
\end{array}
\right),\ \
\partial_{y} = \left(
\begin{array}{ccc}
0 & 0 & 0 \\
0 & 0 & 1 \\
0 & 0 & 0
\end{array}
\right),\ \
\partial_{z} = \left(
\begin{array}{ccc}
-e^{-z} & 0 & 0 \\
0 & e^{z} & 0 \\
0 & 0 & 0
\end{array}
\right),
$$
which are the left translations of the basic vectors:
$$
L_{q*}(e_{1})=e^{-z} \partial_{x}, \ \
L_{q*}(e_{2})=e^{z} \partial_{y}, \ \ L_{q*}(e_{3})
=\partial_{z}.
$$
Since the metric is left-invariant, we have
$$
g_{ij}(x,y,z) = \left(
\begin{array}{ccc}
e^{2z} & 0 & 0 \\
0 & e^{-2z} & 0 \\
0 & 0 & 1
\end{array}
\right).
$$

For the basis  $a_{1}, a_{2}, a_{3}$ we have
$$
L_{q*}(a_{1})=e^{-z}\partial_{x}+e^{z} \partial_{y},
 \ \ L_{q*}(a_{2})=e^{-z}\partial_{x}-e^{z}\partial_{y},
 \ \  L_{q*}(a_{3}) =\partial_{z}.
$$
The inner product takes the form
\beq
\label{e7}
\langle L_{q*}(a_{i}),L_{q*}(a_{j})\rangle =
\langle a_{i},a_{j} \rangle = \left(
\begin{array}{ccc}
2 & 0 & 0 \\
0 & 2 & 0 \\
0 & 0 & 1
\end{array}
\right)
\eeq

In this paper we study the sub-Riemannian problem on the
three-dimensional Lie group $\solv$ defined by the distribution
$\Delta =\mathrm{span} \{a_{1},a_{3}\}$ with metric (\ref{e7}).

Let $G = \solv$,
$\mathcal{G}$ be its Lie algebra with the basic vectors $a_{1}, a_{2},
a_{3}$ (\ref{e3}). We split the Lie algebra $\mathcal{G}$ into the sum
sum $p\bigoplus k$, where $p=\mathrm{span} \{a_{1},a_{3}\},\ \
k=\mathrm{span} \{a_{2}\}$.

Let us consider a two-dimensional left-invariant
distribution $\Delta=\mathrm{span} \{a_{1},a_{3}\}$ in $TG$,
and a left-invariant Riemannian metric
(\ref{e7}) for which the spaces $p$ and $k$ are orthogonal,
 i.e., the metric tensor splits
as follows:
$$
g =\left( g_{ij} \right) =   g_p + g_k.
$$
Let us introduce a parameter $\tau$ and consider the metrics
$$
g_\tau = g_p + \tau g_k.
$$
Every such a metric together with $\Delta$ defines the same
sub-Riemannian manifold because only the restriction
 of the metric onto $\Delta$ is important.

However the Hamiltonian function for the geodesic flows
of these metrics depends on $\tau$:
$$
H(x,p,\tau) = \frac{1}{2}g^{ij}_\tau(x) p_i p_j,
$$
where $g_{ij}g^{jk} = \delta_i^k$.
We have
$$
g_{\tau,ij}=\left(
\begin{array}{ccc}
\frac{1+\tau}{2} e^{2z} & \frac{1-\tau}{2} & 0
\\ \\
\frac{1-\tau}{2}&\frac{1+\tau}{2} e^{-2z} & 0
\\ \\
0 & 0 & 1
\end{array}
\right),
\ \ \
g^{ij}_\tau=\left(
\begin{array}{ccc}
\frac{1+\tau}{2\tau} e^{-2z} & -\frac{1-\tau}{2} & 0
\\ \\
-\frac{1-\tau}{2}&\frac{1+\tau}{2\tau} e^{2z} & 0
\\ \\
0 & 0 & 1
\end{array}
\right).
$$
The Hamiltonian function $H$ for the normal geodesic flow of
the sub-Riemannian metric is obtained from $H(x,p,\tau)$ in the limit
$$
\tau \to \infty,
$$
and we derive
\beq \label{e9}
H(x,y,z,p_{x},p_{y},p_{z})=\frac{1}{4}\,e^{-2z}p_{x}^{2}+
\frac{1}{2}\,p_{x}p_{y}+\frac{1}{4}\,e^{2z}p_{y}^{2}+\frac{1}{2}p_{z}^{2}.
\eeq
The Hamiltonian equations $\dot{x}^i = \frac{\partial H}{\partial p_i},
\dot{p}_i = - \frac{\partial H}{\partial x^i}$
take the form
\beq
\label{e10}
\begin{array}{lll}
\dot{x}=\dfrac{1}{2}\,e^{-2z}p_{x}+\dfrac{1}{2}\,p_{y}, &  \ \  & \dot{p_{x}}=0, \\ & \\
\dot{y}=\dfrac{1}{2}\,e^{2z}p_{y}+\dfrac{1}{2}\,p_{x}, &  \ \  &
\dot{p_{y}}=0,
\\ & \\
\dot{z}=p_{z}, &  \ \  &
\dot{p_{z}}=\dfrac{1}{2}\,e^{-2z}p_{x}^{2}-\dfrac{1}{2}\,e^{2z}p_{y}^{2}.
\\ & \\
\end{array}
\eeq

These differential equations can be derived
from the Pontryagin maximum principle.
The corresponding Hamiltonian takes the form
$$
H(x,y,z,p_{x},p_{y},p_{z},p_{0},u_{1},u_{3})
=\frac{1}{\sqrt{2}}\left(u_{1} p_{x} e^{-z}+u_{1} p_{y} e^{z}\right)+u_{3} p_{z}+p_{0}(u_{1}^{2}+u_{3}^{2}),
$$
where $p_{0}=-\dfrac{1}{2}$, $u_{1}, u_{3}$ are control functions.

The system (\ref{e10}) has three first integrals:
$$
I_{1}=H,\ \ I_{2}=p_{x},\ \ I_{3}=p_{y},
$$
which are functionally independent almost everywhere,
and therefore the system is completely integrable.

Since the flow is left-invariant as well as the distribution
 $\Delta$ and the metric, without loss of generality
we assume, that all geodesics
originate at the unit of group, that is, we have the following
initial conditions for the system (\ref{e10}):
\beq
\label{e11}
x(0)=0,\ \ y(0)=0,\ \ z(0)=0.
\eeq
 In the sequel, we put
$$
H=\frac{1}{2},\ \ \frac{p_{x}}{\sqrt{2}}=a,\ \ \frac{p_{y}}{\sqrt{2}}=b.
$$
By substituting these expressions into (\ref{e9}), we obtain
\beq
\label{e12}
1=\left(e^{-z} a+e^{z} b\right)^{2}+p_{z}^{2},
\eeq
which implies
$$
p_{z}=\pm \sqrt{1-\left(e^{-z}\, a+e^{z}\, b\right)^{2}}.
$$

By substituting this expression to the third
equation of (\ref{e10})
we obtain equation for the temporal variable $t$
for positive values of $p_{z}$
\beq \label{e13}
t=\int\frac{dz}{\sqrt{1-\left(e^{-z}\, a+e^{z}\, b\right)^{2}}}.
\eeq
If $p_{z}<0$, then all calculations will be similar, but with contrary sign.

Let us make the change of variables
$$
u=e^{z},
$$
and rewrite (\ref{e13})  as
\beq
\label{e14}
t=\int \frac{du}{\sqrt{u^{2}-\left( a+bu^{2}\right)^{2}}}.
\eeq
The last expression is not integrated in terms
of elementary functions and defines an
elliptic integral, except of special cases,
when this elliptic integral degenerates.
These cases will be discussed below.

Let us consider the generic case $a\neq0$ and $b\neq0$.

The subradical expression in (\ref{e14}) haves
discriminant $D=1-4ab \geq 0$ accordingly (\ref{e12}).

$D = 0$ if and only if $p_{z} = 0$ accordingly system (\ref{e10})
and equation (\ref{e12}).
That case is degenerative.

Thus, if $D>0 \left(a b<\dfrac{1}{4}\right)$, then there exist $\sigma_{1}^{2}$ and $\sigma_{2}^{2}$,
such that the following
$$
u^{2}-\left( a+bu^{2}\right)^{2}=-b^{2}u^{4}+(1-2ab)u^{2}-a^{2}=
$$
$$
=-b^{2}(u^{2}-\sigma_{1}^{2})(u^{2}-\sigma_{2}^{2})=
\sigma_{1}^{4}b^{2}\left(1-\frac{u^{2}}{\sigma_{1}^{2}}\right)\,
\left(\frac{u^{2}}{\sigma_{1}^{2}}-\frac{\sigma_{2}^{2}}{\sigma_{1}^{2}}\right),
$$
and
\beq
\label{sigma}
\sigma_{1, 2}^2 = \frac{1-2 a b \pm \sqrt{1-4 a b}}{2 b^2}
\eeq

Put
\beq\label{e33}
w=\frac{u}{\sigma_{1}}
\eeq
and rewrite (\ref{e14}) in the following form
\beq\label{e34}
t=\frac{1}{\sigma_{1}b}\, \int
\frac{d w}{\sqrt{\left(1-w^{2}\right)\,
\left(w^{2}-\dfrac{\sigma_{2}^{2}}{\sigma_{1}^{2}}\right)}}.
\eeq
We apply the Jacobi elliptic function (\ref{ell6}) in order to inverse this integral:
$$
\sigma_{1}b\, t=\int_{1}^{\dn (\sigma_{1}bt)}\frac{dw}{\sqrt{\left(1-w^{2}\right)\,
\left(w^{2}-\frac{\sigma_{2}^{2}}{\sigma_{1}^{2}}\right)}},
$$
where $k'^{2}=\frac{\sigma_{2}^{2}}{\sigma_{1}^{2}}$. Therefore
$$
w = \dn(\sigma_{1}bt,k),
$$
where
\beq
\label{kmod}
k^2 = 1 - \frac{\sigma_2^2}{\sigma_1^2}.
\eeq
By inverting (\ref{e33}), putting $u=e^{z}$, and keeping in mind the initial
condition (\ref{e11}) and the equality $\dn(0, k)=1$, we obtain
$$
z(t)=\ln \dn(\sigma_{1}b t,k).
$$
By substituting this expression into the first equation of (\ref{e10})
and integrating it in elliptic functions (see\cite{GR}), we derive:
$$
x(t)=\frac{1}{\sqrt{2}}\left[\frac{a}{\sigma_{1}b}\left(-\frac{k^{2}\sn
(\sigma_{1}bt)\, \cn (\sigma_{1}bt)}{k'^{2}\sqrt{1-k^{2}\sn^{2}
(\sigma_{1}bt)}}+\frac{1}{k'^{2}}\, \E(\am(\sigma_{1}bt),\,
k)\right)+bt\right]+C,
$$
where $\E(x,k)$ is elliptic integral of second kind (\ref{ell2}).

Since $\sn(0,k)=0, \cn(0,k)=1, \am(0,k)=0$ and $\E(0,k)=0$, we have $C=0$.
From the second equation of this system we conclude that
$$
y(t)=\frac{1}{\sqrt{2}}\left(\frac{1}{\sigma_{1}}\, \E(\am(\sigma_{1}bt),\, k)+at\right)+Q,
$$
with $Q = \mathrm{const}$. By (\ref{e11}), we compute that $Q=0$.


Let us now consider the cases, when the elliptic integral (\ref{e14}) degenerates:
\begin{enumerate}
\item
$a=0,\ \ b=0$;

\item
$a=0,\ \ b\neq 0$;

\item
$b=0,\ \ a\neq 0$;

\item
$D=0 \left(a b=\dfrac{1}{4}\right)$.
\end{enumerate}

We consider its successively:

\textbf{1)} $a=0,\ \ b=0$.
From the equations (\ref{e10}), (\ref{e11}) and (\ref{e13}) is cleaner, that
\beq
\label{e15}
x(t)=0, \ \ y(t)=0, \ \ z(t)=t.
\eeq

\textbf{2)} $a=0,\ \ b\neq 0$. We have $p_{x}=0,\ \ p_{y}=\sqrt{2}b$.\\
The equation (\ref{e14}) is rewritten as
$$
t=\int\frac{du}{\sqrt{u^{2}-b^{2}u^{4}}}.
$$
By integration and transformation the resulting expression by
the inverse change of variable, we obtain
$u=e^{z}$:
$$
e^{t}= \frac{Cbe^{z}}{1+\sqrt{1-b^{2}e^{2z}}},
$$
where $C=\mathrm{const}$ and $C>0$.
The last expression together with the initial condition
(\ref{e11}) implies
$$
C=\frac{1+\sqrt{1-b^{2}}}{b},
$$
and we derive that
$$
e^{z}=\frac{2Ce^{t}}{b\left(C^{2}+e^{2t}\right)},
$$
i.e.,
\beq
\label{e16}
z(t)=\ln\frac{2Ce^{t}}{b\left(C^{2}+e^{2t}\right)},
\eeq
which after substituting the formula for $C$ takes the form
$$
z(t)=\ln
\frac{2(1+\sqrt{1-b^{2}})\, e^{t}}
{2(1+\sqrt{1-b^{2}})-b^{2}+b^{2}e^{2t}}.
$$
By the first
equation of (\ref{e10}), we have
$$
x(t)=\frac{b}{\sqrt{2}}\, t,
$$
and the second equation of (\ref{e10}) together with
(\ref{e16}) and (\ref{e11}) implies
$$
y(t)=-\frac{\sqrt{2}C^{2}}{b(C^{2}+e^{2t})}+\frac{\sqrt{2}C^{2}}{b(C^{2}+1)}.
$$
Finally in the case \textbf{2)} we have the explicit formulas for solutions:
\beq
\label{e20}
\begin{array}{l}
x(t)=\dfrac{b}{\sqrt{2}}\, t \\ \\
y(t)=-\dfrac{\sqrt{2}\left(2(1+\sqrt{1-b^{2}})-b^{2}\right)}
{2b(1+\sqrt{1-b^{2}})-b^{3}+b^{3}e^{2t}}+
\dfrac{\sqrt{2}\left(2(1+\sqrt{1-b^{2}})-b^{2}\right)}{2b(1+\sqrt{1-b^{2}})} \\ \\
z(t)=\ln\dfrac{2(1+\sqrt{1-b^{2}})\,
e^{t}}{2(1+\sqrt{1-b^{2}})-b^{2}+b^{2}e^{2t}}.
\end{array}
\eeq

\textbf{3)} $b=0,\ \ a\neq 0$. We have $p_{y}=0,\ \ p_{x}=\sqrt{2}a.$\\
The equation (\ref{e14}) takes the form
$$
t=\int\frac{du}{\sqrt{u^{2}-a^{2}}}.
$$
Let us put $u=e^{z}$ and derive
$$
e^{t}=(e^{z}+\sqrt{e^{2z}-a^{2}})
\, C,
$$
where $C = \mathrm{const}$ and $C>0$. The last expression together with
(\ref{e11}) implies
\beq
\label{e22}
C =\frac{1}{1+\sqrt{1-a^{2}}},
\eeq
from which we obtain
$$
z(t)=\ln\frac{C^{2}a^{2}+e^{2t}}{2Ce^{t}},
$$
where $C$
is given by (\ref{e22}).
As in the case 1) we derive from the first two
equations of (\ref{e10}) that
\beq
\label{e24}
x(t)=-\frac{\sqrt{2}aC^{2}}{e^{2t}+C^{2}a^{2}}+
\frac{\sqrt{2}aC^{2}}{1+C^{2}a^{2}}, \ \
y(t)=\frac{a}{\sqrt{2}}\, t.
\eeq
Finally we obtain
\beq\label{e26}
\begin{array}{l}
x(t)=-\dfrac{\sqrt{2}\,
a}{e^{2t}\left[2(1+\sqrt{1-a^{2}})-a^{2}\right]
+a^{2}}+\dfrac{\sqrt{2}\, a}{2(1+\sqrt{1-a^{2}})}, \\
y(t)=\dfrac{a}{\sqrt{2}}\, t,
\\
\\
z(t)=\ln\left(\dfrac{a^{2}}{2(1+\sqrt{1-a^{2}})\, e^{t}}
+\dfrac{(1+\sqrt{1-a^{2}})\, e^{t}}{2}\right).
\end{array}
\eeq

\textbf{4)} $D=0 \left(a b=\dfrac{1}{4}\right)$.

Note that, the formula (\ref{e12}) by (\ref{e11}) is rewritten as
\beq\
\label{e28}
(a+b)^{2}+p_{z}^{2}=1,
\eeq
which means that
\beq
\label{e29}
|a+b|\leq 1.
\eeq
Then it is clear that  $a=b=\frac{1}{2}$ or $a=b=-\frac{1}{2},$
and for these values the equation (\ref{e28}) implies that $p_{z}=0$.
Therefore solutions to (\ref{e10}) in the case \textbf{4)} are linear:
\begin{equation}
\label{e30}
\begin{split}
x(t)=\frac{t}{\sqrt{2}}, \ \  \ \ y(t)=\frac{t}{\sqrt{2}}, \ \  \ \ z(t)=0, \ \  \ \ a=b=\frac{1}{2};
\\
x(t)=-\frac{t}{\sqrt{2}}, \ \  \ \ y(t)=-\frac{t}{\sqrt{2}},  \ \  \ \ z(t)=0,\ \  \ \ a=b=-\frac{1}{2}.
\end{split}
\end{equation}

Thus we have the following

\begin{theorem}
 In a generic case
the normal geodesics (with the initial condition (\ref{e11}))
are described by the formulas (for $p_{z}>0$):
\beq
\label{e35}
\begin{array}{l}
x(t)=-\dfrac{ak^{2}\sn (\sigma_{1}bt)\, \cn (\sigma_{1}b
t)}{\sqrt{2}\sigma_{1}bk'^{2}\sqrt{1-k^{2}\sn^{2}
(\sigma_{1}bt)}}+\dfrac{a\, \E(\am(\sigma_{1}bt),\, k)}{\sqrt{2}\sigma_{1}bk'^{2}}+\dfrac{b}{\sqrt{2}}t \\
\\
y(t)=\dfrac {\E(\am(\sigma_{1}bt),\, k)}{\sqrt{2}\sigma_{1}}+\dfrac{a}{\sqrt{2}}t \\
\\
z(t)=\ln \dn(\sigma_{1}bt),
\end{array}
\eeq
where the parameters $\sigma_1$ and $k$ are determined
by $a$ and $b$ $(a b<\frac{1}{4})$ via (\ref{sigma}) and (\ref{kmod}).

In the degenerated cases 1--4 the normal geodesics
(with the initial condition (\ref{e11}))
are described in terms of elementary functions
by the formulas (\ref{e15}), (\ref{e20}), (\ref{e26}), and (\ref{e30}).
\end{theorem}

Notice that normal geodesics in the theorem are parameterized by $a, b$.
Constants $k, \sigma_1$ are defined by $a, b$ as explained before.

The qualitative behavior of generic normal geodesic is quite
complicated. Figures 1 and 2 show parts of the geodesic spheres of
radius of 0.15 and 0.25 (a scale on each figure itself; axis $z$ is
exponentially scaled). A grid on the spheres corresponds to two
parameters $\theta$ and $\mu$, where $\theta$ is angle of the
initial vector of the geodesic with respect to the axis $x$ and
$\mu$ is the initial acceleration value $x+y$ along the geodesic,
i.e. $\mu$ can be interpreted as the acceleration with which the
geodesic is drawn out of the starting point. In Figures $\theta$
varies from $\pi/6$ to $5 \pi/6$ (part of the sphere of $ - \pi / 6
$ to $ -5 \pi / 6 $ is obtained as a mirror). Parameter $\mu$ varies
from $-45$ to $45$.

On this grid it can be seen only the qualitative behavior of a sphere with
increasing radius. The drawings practically do not show those parts of spheres,
which are too fast going to infinity, as well as those, which
coincide to the geodesics, changing too quickly the direction.
We can see that part of geodesics
starting at small angle to plane $x, y$  goes to large values of coordinates $x, y$ very quickly,
even for not large values of parameter $\mu$. For sufficiently large $\theta$ and mean values of
$|\mu|$ geodesics deviate not too much from plane $x=y$, but if $|\mu|$ increases the deviation from this plane begins.
Assuming exponential scale of axis $z$ we see that coordinate $z$ increases much more
slow than $x$ and $y$.

\begin{figure}[ht]
\begin{center}
\epsfig{file=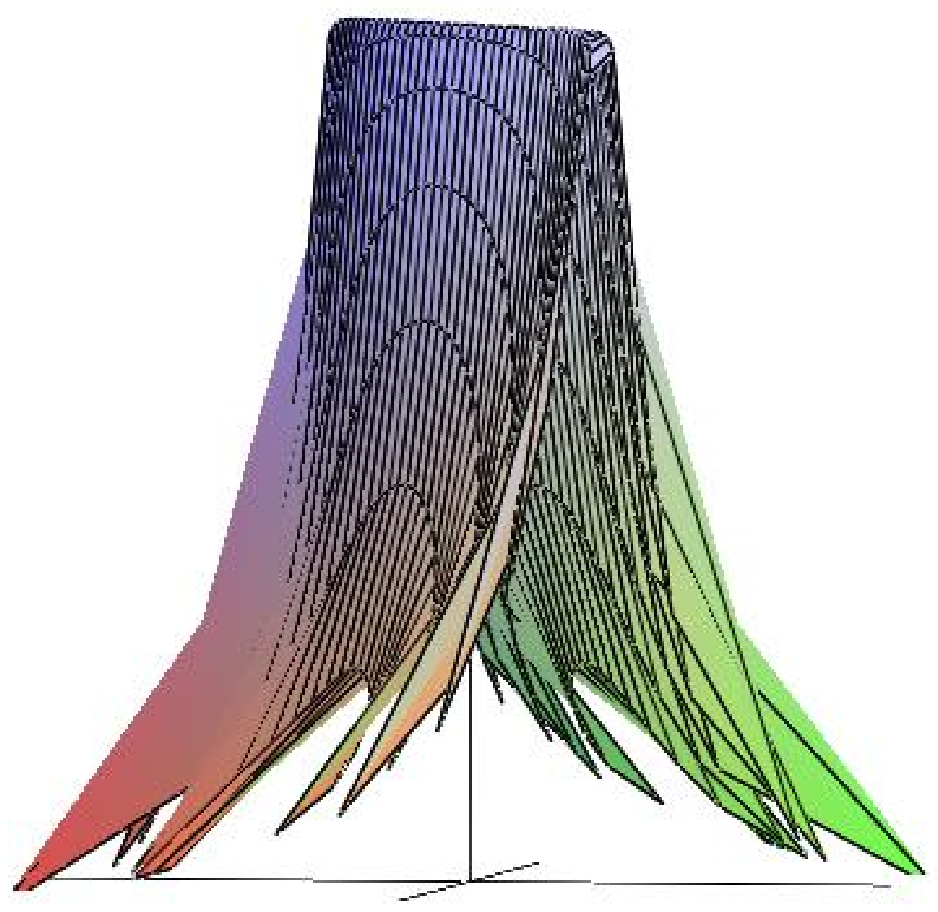,height=80mm,clip=} \hskip1cm
\caption{r=0,15}
\end{center}
\end{figure}

\begin{figure}[ht]
\begin{center}
\epsfig{file=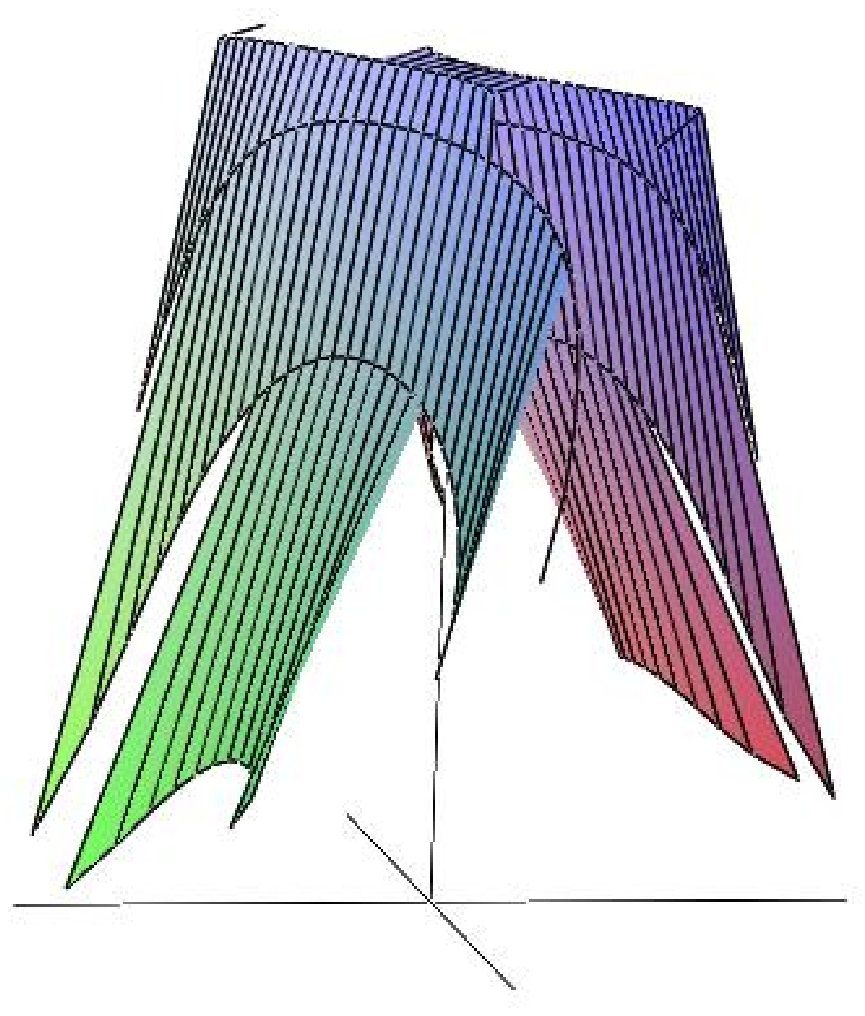,height=80mm,clip=} \hskip1cm \caption{
r=0,25 }
\end{center}
\end{figure}

\end{document}